\documentclass[11pt]{article}
\usepackage{latexsym,euscript,amsmath,amssymb,amsbsy,amsfonts,amsthm,amsopn,amstext,amsxtra,amscd}
\usepackage{epsfig}
\usepackage{graphics}
\usepackage{url}

\usepackage[english]{babel}

\theoremstyle{plain}
\newtheorem{theorem}{Theorem}

\newtheorem{point}{}

\theoremstyle{definition}
\newtheorem{example}{Example}

\newcommand{\MM}{{\mathbb M}}
\newcommand{\N}{{\mathbb N}}

\newcommand{\Z}{{\mathbb Z}}

\newcommand{\D}{\text{$\mathcal{D}$}}

\newcommand{\M}{{\mathbb M}}

\begin{document}
\centerline{\bf METRICS ON $\N$  AND}
\centerline{\bf THE DISTRIBUTION OF SEQUENCES}
\vskip0,5cm
\centerline{Milan Pa\v st\'eka}
\vskip0,5cm
{\bf Abstract.} {\it In the first part of this paper the notion of natural metric on the set of natural numbers is defined. It is such metric that the completion of $\N$  is a compact metric space that a probability borel measure exists in order that the sequence $\{n\}$ is uniformly distributed.
A necessary and sufficient condition where a given metric is
natural. Later we study the properties of sequences uniformly continuous with respect to given natural metric.
Inter alia Theorems \ref{continuity of df} and
\ref{cadf1} characterise the continuity of distribution function.}


\section{Introduction}
This paper is continuation of \cite{P-T}, \cite{P-T1} and \cite{pasclt}.

Hermann Weyl was the first to study the distribution of bounded sequences of real numbers in his famous work \cite{WEY} published in 1916. The crucial result of this paper is now known as "Weyl criterion" generalized later for many cases. In 1928 this concept was essentially generalized by Schoenberg in the paper \cite{Sch}. Schoenberg define the asymptotic distribution function of given sequence $\{v(n)\}$ as the limit
$$
F(x) = \lim_{N \to \infty} \frac{|\{n \le N; v(v)< x\}|}{N}
$$
if the limit exists for each $x$ - real number. The notion of asymptotic density allows define the asymptotic distribution function analogously as the distribution function of random variable in probability space. This leads to the possibility of application of results of Edmund Hlawka, \cite{HLA}, and other authors which  defined later the uniform distribution for compact metric space with Borel probability measure.
From the 70 ies the monographs \cite{Hla}, \cite{K-N}, \cite{D-T}, \cite{SP} and \cite{Str} were published. The aim of this paper is to apply such a probability  measure for describe the distribution of sequences in some cases using asymptotic density. The properties of the compact ring of polyadic integers were utilized in  \cite{P-T1}. The additive compact group of polyadic integers provides the existence of Haar probability measure. This measure plays an important rule in the study of distribution of polyadic continuous sequences. Very useful is the fact that
the sequence $\{n\}$ is uniformly distributed in the probability space of polyadic integers in the sense of Hlawka.
In generally the completion of $\N$ is not necessary an algebraic structure. The other name for asymptotic density is "natural density". Because of this we call the corresponding metric and measure "natural".
We define the "natural metric" as metric that the completion is a compact providing the existence of such probability measure that the mentioned sequence of natural numbers in natural order is uniformly distributed. Such measure is called "natural measure". By the application of Riesz representation theorem to the linear space of continuous functions we derive a necessary and sufficient condition providing that given metric is natural (\ref{natural} and Theorem \ref{natural1} ).

The connection between the asymptotic distribution function and the distribution function of random variable is proven in
\cite{P-T1} proven via Buck's measure density (Theorem 6 in \cite{P-T1}). The transfer for general situation allows the
"measure density"(see \cite{PAS5}). This connection is described in Theorem \ref{equivalenceI}. These statements then can by applied to statistical independence of sequences (introduced and studied by Gerard Rauzy in his monograph \cite{Ra})and study the independence of sequences as independence of random variables.  \newline
{\bf Notation.} We shall use the following notation: \newline $\N$ - the set of natural numbers \newline $\mathcal{C}(\mathfrak{d})$ - the linear space of sequences uniformly continuous with respect to given metric $\mathfrak{d}$ equipped with supremum norm \newline $\mathcal{L}(M)$ - linear subspace generated by the elements of the sets $M$ \newline
$\mathcal{X}_S$ - indicator function of the set $S$ \newline
$cl(S)$ - closure of the set $S$ in given metric space \newline
$int(S)$ - interior of the set $S$ in given metric space
\newline
$\partial(S)= cl(S) \setminus int(S)$
\newline
$\mathcal{B}(M)$ the system of all Borel subsets of given metric space $M$.
\newline
$\mathcal{B}_c(M)$ the sets of all sets of $P$ - continuity
with respect to given probability measure $P$.
 \newline $E_N(v)=\frac{1}{N}\sum_{n=1}^N v(n)$ \newline
$E(v)=\lim_{N \to \infty} E_N(v)$ (if it exists)\newline
$D_N^2(v)= E_N((v-E(v))^2)$ \newline
$D^2(v)=\lim_{N\to \infty} D^2_N(v)$ (if it exists)\newline
$\underline{d}(S) \liminf_{N\to \infty} E_N(\mathcal{X}_S)$ - lower asymptotic density of the set $S$
\newline
$\overline{d}(S)=\limsup_{N\to \infty} E_N(\mathcal{X}_S)$ - upper asymptotic density of the set $S$
\newline
$d(S)= E(\mathcal{X}_S)$ - asymptotic density of the set $S$ (if it exists)
\newline
$\D$ - the system of all subsets of $\N$ having asymptotic density \newline
$\nu^\ast(S)=P(cl(S))$
\newline
$\D_\nu$ - the system of all $\nu^\ast$ - measurable sets \newline
$\nu$ - the restriction of $\nu^\ast$ on $\D_\nu$
\newline
$f(v)=\{f(v(n))\}$ for the sequence $v=\{v(n)\}$ \newline
$|I|$ - length of interval $I$ \newline
$\lambda(S)$ - Lebesque measure of $S$, (for interval $I$ we have $\lambda(I)=|I|$) \newline
$\lambda_F(S)$ - Lebesque - Stjelties measure of $S$ given by the function $S$

\section{Preliminaries}
We start by recalling the basic concepts used in whole paper.

{\bf Asymptotic density.} Let $A \subset \N$. The value
$$
\limsup_{N \to \infty} E_N(\mathcal{X}_A)= \limsup_{N \to \infty} \frac{|A \cap [1,N]|}{N}= \overline{d}(A)
$$
is called the {\it upper asymptotic density} of the set $A$. It is easy to check that:
\begin{equation}\overline{d}(\N)=1, \overline{d}(S)=0
\end{equation} if $S$ is a finite set, and
\begin{equation}
\label{sub}
\overline{d}(S_1\cup S_2)\le \overline{d}(S_1) + \overline{d}(S_2)
\end{equation}
for $S_1, S_2 \subset \N$.

A system of subsets of $\N$ denoted by $\mathcal{Y}$ is called
the {\it q - algebra} if and only if : \newline
i) $\N \in \mathcal{Y}$, \newline
ii) $S \in \mathcal{Y} \Rightarrow \N \setminus A \in \mathcal{Y}$, \newline
iii) $S_1, S_2 \in \mathcal{Y} \land S_1 \cap S_2 = \emptyset
\Rightarrow S_1 \cup S_2 \in \mathcal{Y}$.

A finitely additive probability measure $\pi$ defined on
q - algebra $\mathcal{Y}$ is called {\it density} if and only if \newline
iv) $S \in \mathcal{Y} \Leftrightarrow \forall \varepsilon >0
\exists S_1, S_2 \in \mathcal{Y}; S_1 \subset S \subset S_2 \land \pi(S_2) - \pi(S_1) < \varepsilon.$

Let $\D$ be the system of all sets $A \subset \N$ such that
$$
\overline{d}(A)+\overline{d}(\N \setminus A)=1.
$$
\begin{point} A set $A$ belongs to $\D$ if and only if
the limit $\lim_{N \to \infty} E_N(\mathcal{X}_A) := d(A)$
exists. $\D$ is a q - algebra and $d$ is density on $\D$.
\end{point}
The value $d(A)$ is called the {\it asymptotic density} of $A$, and in this case we say also that $A$ has asymptotic density. From (\ref{sub}) we get immediately
\begin{point}
\label{sub1} A set $A \subset \N$ belongs to $\D$ if and only if $\overline{d}(A)+\overline{d}(\N \setminus A) \le 1$.
\end{point}

{\bf Distribution function of sequence.}
Schoenberg's definition in \cite{Sch} can be rewritten as follows: \newline
a bounded sequence $\{v(n)\}$ has {\it asymptotic distribution function} if and only for each real number $x$ the set $v^{-1}([0,x))$ belongs to $\D$. In this case the function
\begin{equation}
\label{asdf}
F(x)= d(v^{-1}([0,x)))
\end{equation}
is called the {\it asymptotic distribution function} of
$\{v(n)\}$.
This concept was later studied by a number of authors.

The following criterion is well known as {\it Weyl criterion}:
\begin{theorem}
\label{wk} Let $\{v(n)\}$ be a sequence of elements of closed interval $[a,b],$\newline $a<b$. Suppose that $F$ is continuous non decreasing function such that $F(a)=0, F(b)=1$. Then $F$ is the asymptotic distribution function of $\{v(n)\}$ if and if
\begin{equation}
\label{WKI}
\lim_{N\to \infty} E_N(f(v))= \int_a^b f(x)dF(x)
\end{equation}
for each continuous function $f$ defined on $[a,b]$.
\end{theorem}
For the proof we refer to \cite{D-T}, \cite{K-N}, \cite{Str} .

A sequence $\{v(n)\}, v(n) \in [0,1]$ is called {\it uniformly
distributed modulo} $1$ if and only if it has asymptotic distribution function $F$, where $F(x)=0$ for $x \le 0$,
$F(x)=x$ for $x \in [0,1]$ and $F(x)=1$ for $x\ge 1$, (see \cite{WEY}).
Since each continuous function defined on $[0,1]$ can be uniformly approximated by polynomial the following form
of Weyl criterion can be proven (see \cite{D-T}, \cite{K-N}, \cite{Str}):
\begin{point}
\label{wkpolynom} A sequence $\{v(n)\}, v(n) \in [0,1]$ is uniformly distributed modulo $1$ if and only if for each
$h \in \N$ we have
$$
\lim_{N \to \infty} E_N(v^h)=\frac{1}{h+1}.
$$
\end{point}

There are lot of functions $F$ with points of discontinuity.
For this reason we will use a slightly more general definition.

We say that $F$ is {\it weak asymptotic distribution functon}
of $\{v(n)\}$ if and only if the equality (\ref{asdf}) holds for each $x \in [a,b]$ - point of continuity of $F$.

The complete characterization of sequences having weak asymptotic distribution function is made in \cite{Str} page 33.\footnote{ A bounded sequence $\{v(n)\}$ of elements of $[0,1]$ has weak asymptotic distribution function if and only if
$$
\lim_{M,N \to \infty} \frac{1}{MN}\sum_{n\le N, m\le M}|v(n)-v(m)|-
$$
$$
-\frac{1}{2N^2}\sum_{n\le N} |v(n)-v(m)|-\frac{1}{2M^2}\sum_{m\le M}|v(n)-v(m)|=0.
$$
By the linear transformation can be this result extended for arbitrary interval $[a,b]$.
}

{\bf Statistical independence.}
A sequences $\{v_j(n)\}, j=1,\dots, k$ are called
{\it statistical independent} if and only if for every
continuous functions $f_1, \dots, f_k$ we have
$$
E_N(f_1(v_1)\dots f_k(v_k))- \prod_{j=1}^k E_N(f_j(v_j)) \to 0
$$
for $N \to \infty$.

Let us consider a sequence of vectors $\{(v_1(n), \dots, v_k(n))\}$. We say that this sequence has {\it asymptotic distribution function} $F(x_1, \dots, x_k)$ if and only if
for each $x_1, \dots, x_k$ the set $\cap_{j=1}^k v_j^{-1}((-\infty, x_j))$ belgs to $\D$ and its asymptotic density is equal $F(x_1, \dots, x_k)$.

\begin{point}
\label{INDEP} If the sequences $\{v_1(n)\}, \dots, \{v_k(n)\}$ of elements of some interval $[a,b],$ \newline $a<b$
have continuous asymptotic distribution function then these
sequences are statisticaly independent if and only if
$$
F(x_1, \dots, x_k)= F_1(x_1)\dots F_k(x_k)
$$
for $x_1, \dots, x_k$ real numbers, where $F_j$ is asymptotic distribution function of the sequence $\{v_j(n)\}, j=1, \dots, k$.
\end{point}
(Theorem 3 of \cite{P-T1}.)

\begin{point} If the sequences $\{v_1(n)\}, \dots, \{v_k(n)\}$ satisfy the conditions of \ref{INDEP} then the sequence
$\{v_1(n)+\dots + v_k(n)\}$ has continuous asymptotic distribution function.
\end{point}

{\bf Uniform distribution in compact space.}
Let $(\M, \mathfrak{d}, P)$ be a compact metric space. A sequence $\{v(n)\}, v(n) \in \M$ is {\it uniformly distributed in $\M$} if and only if
\begin{equation}
\lim_{N\to \infty} E_N(f(v))=\int f dP
\end{equation}
for each continuous function $f$ defined on $\M$ (see \cite{HLA}).

Let us remark that important argument in this topics is the regularity of measure. A borel probability measure $P$ on topological space is called {\it regular} if and only if
for each borel set $S$ we have
$$
P(S)=\inf\{P(O); S \subset O, O - \text{open set}\}=
\sup\{P(O); C \subset S, O - \text{closed set}\}.
$$
Following result is proven in \cite{par}:
\begin{point}
\label{regularmeasure} Each borel probability measure on metric space
is regular.
\end{point}

{\bf The sets of $P$ - continuity.}
Uniform distribution in compact spaces can be characterized by asymptotic density.

If  $S \subset \M$ then the set $\partial(S):= cl(S)\setminus Int(S)= cl(S)\cap cl(\M \setminus S)$ is called the {\it border} os $S$. A borel set
$S \subset M$ will be called {\it set of $P$ continuity} if and only if $P(\partial(S))=0$.
These sets play similar rule as intervals in the case of uniform distribution mod $1$.

\begin{theorem} A sequence $\{v(n)\}, v(n) \in \M$ is  uniformly distributed in $\M$ if and only if for set of $P$ continuity
$S$ we have $v^{-1}(S) \in \D$ and $d(v^{-1}(S))=P(S)$.
\end{theorem}
For the proof of this result we refer to \cite{K-N}.

A q - algebra $\mathcal{Y}$ is called the algebra of sets if and only if for $S_1, S_2 \in \mathcal{Y}$ we have
$S_1 \cup S_2 \in \mathcal{Y}$, (the condition of disjointness
in iii) is omitted).

\begin{point}
\label{pcontalgebra}
The system of the sets of $P$ - continuity
$\mathcal{B}_c(M)$ is an algebra of sets.
\end{point}

{\bf Proof.} Clearly $\emptyset, M \in \mathcal{B}_c(M)$. Moreover $\partial(M \setminus S)=\partial(S)$ and so
$S \in \mathcal{B}_c(M) \Rightarrow M \setminus S \in \mathcal{B}_c(M)$.

Consider $S_1, S_2 \subset M$.  We have $\partial(S_1 \cup S_2) \subset \partial(S_1) \cup \partial(S_2)$. Thus
$$
S_1, S_2 \in \mathcal{B}_c(M) \Longrightarrow S_1 \cup S_2
\in \mathcal{B}_c(M).
$$
\qed

Denote
$$ B(\alpha, \varepsilon) = \{\beta \in \MM; \mathfrak{d}(\alpha,
\beta) < \varepsilon\}
$$
the open ball with center $\alpha \in \MM$ and radius
$\varepsilon >0$ and
$$ B[\alpha, \varepsilon] = \{\beta \in \MM; \mathfrak{d}(\alpha,
\beta) \le \varepsilon\}
$$
the closure of the open ball with center $\alpha \in \MM$ and radius
$\varepsilon >0$.

\begin{point}
\label{pcont}
For each $\alpha \in \MM$ and $r>0$ the set $B(\alpha, r')$ belongs to $\mathcal{B}_c(M)$  for $r'\in [0,r]$ with except of countable subset of $[0,r]$.
\end{point}

{\bf Proof.} For $\alpha \in M$, $r'>0$ we have
$\partial(B(\alpha, r'))=\{\beta \in \M; \mathfrak{d}(\alpha, \beta)=r'\}$. This yields that for $r'' \neq r'$ the sets
$\partial(B(\alpha, r')), \partial(B(\alpha, r''))$ are disjoint. Each system of disjoint sets with positive measure
is countable and the assertion follows. \qed

From \ref{pcont} we get:
\begin{point}
\label{unionofpc} For each $\alpha \in \M$ and $r >0$ such increasing sequence of $r_n < r$ that $r_n \to r, n \to \infty$ that the open balls $B(\alpha, r_n)$ are the set of
$P$ - continuity.
\end{point}
Since in compact metric exists an open base formed from open balls, we can deduce from \ref{unionofpc} :
\begin{point}\label{openbase}
An open base of the metric space $(\M, \mathfrak{d})$
exists that it contains only the open balls which are the sets
of $P$ - continuity.
\end{point}

\section{Measure density, extension of sequences} Suppose that $(\M,\mathfrak{d})$ is a compact metric space with Borel probability measure $P$, containing $\N$ as dense subset.
For $S \subset \N$ we can define
\begin{equation}
\label{measuredensity}
\nu^\ast(S)=P(cl(S))
\end{equation}
where $cl(\cdot)$ denote topological closure of given set in $\MM$.
Let us denote by
$$
\D_\nu=\{S \subset \N; \nu^\ast(S)+\nu^\ast(\N \setminus S)=1\}
$$
The elements of $\D_\nu$ will be called $\nu^\ast$ - measurable sets.
\begin{example}If $\M$ is the ring of polyadic integers, then
the corresponding measure density is Buck's measure density,
(see \cite{BUC}, \cite{PAS4}, \cite{PAS5}).
\end{example}
\begin{point}
The system $\D_\nu$ is an algebra of sets and $\nu = \nu^\ast
|_{\D_\nu}$ is density on $\D_\nu$.
\end{point}
\begin{point}
\label{measrpcont} A set $S \subset \N$ is a $\nu^\ast$ - measurable set if and only if its closure $cl(S)$ is the set of $P$ - continuity.
\end{point}

{\bf Proof.} Clearly $\MM = cl(S \cup (\N \setminus S))=
cl(S)\cup cl(\N \setminus S)$. Thus
$$
1=P(cl(S))+P(cl(\N \setminus S))-P(cl(S)\cap cl(\N \setminus S)).
$$
And so
$$
1=\nu^\ast(S)+\nu^\ast(\N \setminus S)-P(\partial cl(S)),
$$
and the assertion follows.

To each real valued sequence $\{v(n)\}$ uniformly continuous with respect to $\mathfrak{d}$  we can associate a continuous real function $\tilde{v}$ defined
\begin{equation}
\label{exttorv}
\tilde{v}(\alpha)=\lim_{k \to \infty} v(n_k)
\end{equation}
where $\alpha \in \MM$ and $n_k \to \alpha$ with respect to the extension of the metric $\mathfrak{d}$. The function $\tilde{v}$ is continuous and so is measurable. This allows
to consider it as random variable.

Applying (\ref{exttorv}) and the uniform continuity of $\{v(n)\}$ we get the inclusion
$$
cl((v^{-1}((-\infty, x))) \subset \tilde{v}^{-1}((-\infty, x]))
$$
follows. This yields
\begin{equation}
\label{prvanerovnost}
\nu^\ast(v^{-1}((-\infty, x))) \le P(\tilde{v} \le x).
\end{equation}
From the other side we have for each $x'>x$ we have
$$
\tilde{v}^{-1}((-\infty, x])) \subset cl(v^{-1}((-\infty, x'))),
$$
and so we get
\begin{equation}
\label{druhanerovnost}
P(\tilde{v} \le x) \le \nu^\ast(v^{-1}((-\infty, x'))).
\end{equation}
Analogously we have
$$
cl(v^{-1}([x, \infty)))\subset
\tilde{v}^{-1}([x, \infty)),
$$
this implies
\begin{equation}
\label{tretianerovnost}
\nu^\ast(v^{-1}([x, \infty))) \le P(\tilde{v} \ge x).
\end{equation}
Finally
$$
\tilde{v}^{-1}([x, \infty)) \subset cl(v^{-1}((x'', \infty))
$$
for each $x'' < x$. And so
\begin{equation}
\label{stvrtanerovnost}
P(\tilde{v}\ge x) \le \nu^\ast(v^{-1}((x'', \infty)).
\end{equation}
Let $F(x)=P(\tilde{v} < x)$ be the distribution function of
the random variable $\tilde{v}$.
\begin{theorem}
\label{pcont0} For each point of continuity of $F$ we have
$v^{-1}((-\infty, x))$ belongs to $\D_\nu$ and
$\nu(v^{-1}((-\infty, x)))=F(x)$.
\end{theorem}

{\bf Proof.} We have $F(x)=P(\tilde{v} \le x)$ for each $x$ - point of continuity of $F$. And so from (\ref{prvanerovnost}) we get
$$
\nu^\ast(v^{-1}((-\infty, x))) \le F(x).
$$
From the other hand side we have
$$
\N \setminus v^{-1}((-\infty, x)))= v^{-1}([x, \infty)).
$$
And we from (\ref{tretianerovnost}) we get
$$
\nu^\ast(\N \setminus v^{-1}((-\infty, x)))\le 1-F(x).
$$
We proven
$$
\nu^\ast(v^{-1}((-\infty, x))+\nu^\ast(\N \setminus v^{-1}((-\infty, x)) \le 1.
$$ \qed

We say that a sequence $\{v(n)\}$ has {\it weak} $\nu$ -
{\it distribution function} if and only if such non decreasing
function $G$ exists that the set $v^{-1}((-\infty,x))$
belongs to $\D_\nu$ and $ \nu(v^{-1}((-\infty,x)))=G(x)$ for each $x$ - the point of continuity of $F_0$.
From Theorem \ref{pcont0} we get immediately
\begin{point}
\label{pcont1}
Each sequence uniformly continuous with respect to $\mathfrak{d}$ has weak $\nu$ - distribution function.
\end{point}

Let $F_0(x)=\nu^\ast(v^{-1}((-\infty,x))$ for $x$ real number.
\begin{point}
\label{nucont} For each $x$ - point of continuity of $F_0$ we have $F_0(x)= P(\tilde{v}\le x)$.
\end{point}

{\bf Proof.} Let $x$ be a point of continuity of $F_0$. From (\ref{druhanerovnost}) we get
$$
P(\tilde{v} \le x) \le F_0(x+\varepsilon)
$$
for each $\varepsilon>0$. Taking account the continuity of
$F_0$ in $x$ we get for $\varepsilon \to 0^+$ that
$P(\tilde{v} \le x) \le F_0(x)$.
From (\ref{prvanerovnost}) we see that
$$
P(\tilde{v} \le x) \ge F_0(x).
$$
We proven $F_0(x)=P(\tilde{v} \le x)$.  \qed

We say that the sequence $\{v(n)\}$ {\it has} $\nu$ {\it distribution function} if and only if $v^{-1}((-\infty, x))$ belongs to
$\D_\nu$ for each real number $x$. In this case the function
$\nu(v^{-1}((-\infty, x)))$ is called the $\nu$ {\it distribution function} of $\{v(n)\}$.

\begin{theorem} Let $\{v(n)\}$ be a sequence uniformly continuous with respect to $\mathfrak{d}$. Then the following statements are equivalent: \newline
a) the distribution function of the random variable $\tilde{v}$ is continuous \newline
b) the sequence $\{v(n)\}$ has continuous $\nu$ distribution function. \newline
In this case  the $\nu$ - distribution function and the distribution function of $\tilde{v}$ coincide.
\end{theorem}
{\bf Proof.} The implication a) $\Rightarrow$ b) follows from
Theorem \ref{pcont}. If b) holds and $F_0$ is contiunuous $\nu$ - distribution function of $\{v(n)\}$ then from \ref{nucont} we get $F_0(x)=P(\tilde{v} \le x)$,$x$ - real number. The continuity of $F_0$ yields
$F_0(x)= P(\tilde{v}<x)$. \qed

The above ideas lead also to:
\begin{theorem}
\label{continuity of df} If $\{v(n)\}$ is a bounded sequence uniformly continuous with respect to $\mathfrak{d}$ then it has a continuous $\nu$ - distribution function of and only if
$$
P(\tilde{v}=x)=0
$$
for each real number $x$.
\end{theorem}
\begin{point}
\label{borelset}
If $\{v(n)\}$ is a bounded sequence uniformly continuous with respect to $\mathfrak{d}$ having a continuous $\nu$ - distribution function $F$ then for every $a<b$ we have
$$
P(\tilde{v} \in (a,b))= F(b)-F(a).
$$
\end{point}
This yields by standard method
\begin{theorem}Under assumptions of \ref{borelset} we have borel set on real line
$$
P(\tilde{v} \in S)=\lambda_F(S).
$$
\end{theorem}

\begin{point}
\label{abscont} Let $\{v(n)\}$ be a sequence of elements of
some interval $[a,b], a<b$ having absolutely continuous
$\nu$ - distribution function $F$.
If $f:[a,b] \to [a,b]$ is such continuous function that for each $x \in [a,b]$ we have $\lambda(f^{-1}(\{x\}))=0$, then the sequence $\{f(v(n))\}$ continuous $\nu$ - distribution function.
\end{point}

{\bf Proof.} Clearly $\tilde{f(v)}=f(\tilde{v})$. Thus  for each $x$ we have
$$
P( \tilde{f(v)}=x)= P(\tilde{v} \in f^{-1}(\{x\})=
\lambda_F (f^{-1}(\{x\}))).
$$
The condition of absolute continuity of $F$ implies
$\lambda_F (f^{-1}(\{x\})))=0$. \qed

{\bf Independence.} If $\{v_1(n)\}, \dots, \{v_k(n)\}$ are sequences having $\nu$ distribution function then the fact that $\D_\nu$ is an algebra of sets guaranties that the set
$\cap_{i=1}^k v_i^{-1}((-\infty, x_i))$ belongs to $\D_\nu$. We say that these sequences are $\nu$ - {\it independent} if and only if
\begin{equation}
\label{nindependence}
\nu\Big(\bigcap_{i=1}^k v_i^{-1}((-\infty, x_i))\Big)=\prod_{i=1}^k \nu(v_i^{-1}((-\infty, x_i))).
\end{equation}
\begin{point} Let $\{v_1(n)\}, \dots, \{v_k(n)\}$ be uniformly
continuous sequences with respect to $\mathfrak{d}$, having
continuous $\nu$ - distribution functions. If these sequences are $\nu$ - independent then the random variables
$\tilde{v}_1, \dots, \tilde{v}_k$ are independent.
\end{point}

{\bf Proof.} Let $F_i$ be the $\nu$ - distribution functions of $\{v_i(n)\}, i=1, \dots, k$. Let us consider the real numbers $x_1, \dots, x_k$. Then for each $\varepsilon >0$ we have
$$
cl\Big(\bigcap_{i=1}^k v_i^{-1}((-\infty, x_i -\varepsilon))\Big)
\subset \bigcap_{i=1}^k \tilde{v}_i^{-1}((-\infty, x_i)) \subset cl\Big(\bigcap_{i=1}^k v_i^{-1}((-\infty, x_i +\varepsilon))\Big).
$$
Therefore
$$
\nu\Big(\bigcap_{i=1}^k v_i^{-1}((-\infty, x_i -\varepsilon))\Big) \le P(\tilde{v}_i < x_i, i\le k)\le \nu\Big(\bigcap_{i=1}^k v_i^{-1}((-\infty, x_i +\varepsilon))\Big),
$$
and so applying the $\nu$ - independence of $\{v_i(n)\}$ we get
$$
\prod_{i=1}^k F_i(x_i-\varepsilon) \le P(\tilde{v}_i < x_i, i\le k) \le \prod_{i=1}^k F_i(x_i+\varepsilon).
$$
For $\varepsilon \to 0^+$ we get $P(\tilde{v}_i < x_i, i\le k)
= \prod_{i=1}^k F_i(x_i)$. \qed

\section{Natural metric and continuity}

Let $\mathfrak{d}$ be such metric on the set of natural numbers
$\N= \{1,2,3, \dots, \}$ that the metric space $(\N,\mathfrak{d})$ is totaly bounded.
Denote by $\MM$  its completion.  A  borel probability measure $P$
defined on $\MM$ is called {\it natural} if and only if the following conditions hold: \newline
a)
the sequence $\{n\}$ of all positive integers is $P$ uniformly
distributed in $(\MM, P)$ and \newline b)
 $P(A)>0$ for every non empty open set $A$.

 We say that the metric $\mathfrak{d}$ is {\it natural} if and only in an natural measure on $\MM$ exists. In the following text we shall suppose that $\mathfrak{d}$ is natural metric,
 $\M$ is completion of the metric space $(\M,\mathfrak{d})$ and $P$ is natural measure on $\M$, and $\nu^\ast$ is corresponding measure density.
 \begin{example}One of example of natural metric is the polyadic metric, (see \cite{N}, \cite{N1}, \cite{PAS5}).
 \end{example}
 \begin{example} A more general type of natural metric is constructed in \cite{pasclt}. This construction will be repeated in this paper later - in the proof of \ref{naturalmetric1}.
 \end{example}

 \begin{point}
 \label{ineq} For each $S \subset \N$ the inequality
 $\overline{d}(S) \le \nu^\ast(S)$ holds.
 \end{point}

 {\bf Proof.} The set $cl(S)$ is closed thus from each its cover by open balls exists finite subcover. By \ref{openbase}
 we can suppose that these balls are the sets of $P$ - continuity. By \ref{pcontalgebra} we get that this union is the set of $P$ - continuity also. Thus from regularity
  of $P$ (see \ref{regularmeasure}) we get that for $\varepsilon >0$ an open set $O$ of $P$ - continuity exists that
 $$
 cl(S) \subset O ,\ P(O) < P(cl(S))+\varepsilon.
 $$
 Since $S \subset cl(S)$ we get
 $$
 \overline{d}(S) \le \nu^\ast(S)+\varepsilon.
 $$
 For $\varepsilon \to 0^+$ we get the assertion. \qed

 This implies immediately

\begin{point}
\label{asympdens} If $S$ is $\nu^\ast$ - measurable set then
$S \in \D$ and $\nu(S)=d(S)$.
\end{point}

We say that a sequence is {\it naturally continuous} if and only if such natural metric $\mathfrak{d}$ exists that this sequence is uniformly continuous with respect to $\mathfrak{d}$.
\begin{point}
\label{wadf} Each naturally continuous sequence has weak asymptotic distribution function. If this sequence is dense in some interval then
this weak asymptotic distribution function is increasing on this interval.
\end{point}

{\bf Proof.} If $\{v(n)\}$ is mentioned sequence, then
\ref{pcont} yields that
$$
F(x)=P(\tilde{v}<x)
$$
is weak $\nu$ - distribution function. Thus from \ref{asympdens} this function is weak asymptotic distribution
function of $\{v(x)\}$. Suppose that $x_1 < x_2$ and given sequence is dense in $[x_1,x_2]$. Then such point of continuity $x' < x''$  of $F$ exist that
$x_1 < x' < x'' < x_2$. The set $\tilde{v}^{-1}((x',x''))$
is non - empty open, and so
$F(x'')-F(x') \ge P( x' < \tilde{v} < x'')>0$. Thus
$F(x'')>F(x')$ therefore $F(x_2)>F(x_1)$. \qed

From \ref{asympdens} we get immediately that each $\nu$ - measurable sequence has asymptotic distribution function. The
connection between $\nu$ - distribution function and asymptotic distribution function for a sequence uniformly continuous sequence is described in the following statement:

\begin{theorem}
\label{equivalenceI}
Let $\{v(n)\}$ be a naturally continuous sequence. Let $F$ be continuous function.
Then the following statements are equivalent \newline
a) $\{v(n)\}$ has continuous asymptotic distribution function
$F$ \newline\
b) the distribution of random variable $\tilde{v}$ is $F$.
\end{theorem}

\begin{point}
\label{densint}
If $\{v(n)\}$ is naturally a continuous sequence of interval $[a,b]$, $a< b$, with continuous asymptotic distribution $F$, then for
subinterval $[c,d] \subset [a,b]$ the set
$v^{-1}([c,d])$ belongs to $\D_\nu$ and
$$
d(v^{-1}([c,d]))=\nu(v^{-1}([c,d]))=P(\tilde{v}\in [c,d]).
$$
\end{point}
From Theorem \ref{continuity of df} we get:
\begin{theorem}
\label{cadf1} A naturally continuous bounded sequence
$\{v(n)\}$
has continuous asymptotic distribution function if and only if
$$
P(\tilde{v}=x)=0
$$
for each real number $x$, where $P$ is corresponding natural
measure.
\end{theorem}

{\bf Mean value, dispersion.}
To each sequence of real numbers $\{v(n)\}$ uniformly continuous with respect to $\mathfrak{d}$ we  associate  the {\it mean value}
$$
E(v)=\lim_{N \to \infty}E_N(v) = \int \tilde{v}dP.
$$
where $P$ is an natural probability measure defined on $\MM$

Some results from \cite{P-T} and \cite{P-T1} will be investigated.
From the fact that the measure of each non empty open set is positive we get that each subset $S \subset \MM$ such that
$P(S)=1$ is dense in $\MM$. This implies
\begin{point}
\label{mis1} If a sequence $\{v(n)\}$ is naturally continuous  and $\tilde{v}(\alpha)=0$ for $\alpha \in S$ where $P(S)=1$ then $v(n)=0$ for each $n \in \N$.
\end{point}
\begin{point}
\label{mis2}
If $\{v(n)\}$ is a sequence is naturally continuous, and for each $\varepsilon >0$ we have
$P(|\tilde{v}|> \varepsilon)=0$
then $\tilde{v}(\alpha)=0$ for each $\alpha \in \MM$.
\end{point}

Analogously as the mean value we define {\it dispersion}
$$
D^2(v)= E((v-E(v))^2),
$$
where $\{v(n)\}$ is a sequence uniformly continuous with respect the metric $\mathfrak{d}$.
This lead to the Chebyshev inequality:
$$
P(|\tilde{v}-E(v)|\ge \varepsilon) \le \frac{D^2(v)}{\varepsilon^2}.
$$
And so \ref{mis2} implies
\begin{point} If $\{v(n)\}$ is naturally continuous then $D^2(v)=0$ if and only if this sequence is constant.
\end{point}
Let $\mathcal{S}$ be a family of sequences. We say that given natural metric $\mathfrak{d}$ is {\it common natural metric} for $\mathcal{S}$ if and only if each sequence from $\mathcal{S}$ is uniformly continuous with respect to $\mathfrak{d}$. In this case we say that these sequences are
{\it jointly  naturally continuous}.

Let $\{v(n)\}$ and $\{u(n)\}$ be common naturally continuous
sequences. Suppose that both are non constant. The value
$$
\rho(v,u)=\frac{|E(uv)-E(u)E(v)|}{D(u)D(v)}.
$$
we shall call {\it correlation coefficient} of $\{v(n)\}$ and $\{u(n)\}$.

Denote
$$
a=\frac{E(uv)-E(u)E(v)}{D^2(v)}, b= E(v)-aE(u).
$$
We have
$$
D^2(u(n)-av(n)+b)=(1-\rho(u,v))D^2(u).
$$
And so we can conclude
\begin{point} If the sequences $\{u(n)\}, \{v(n)\}$ are jointly naturally  continuous then
$\rho(v,u)=1$ if and only $u(n)=av(n)+b$ for $n \in \N$.
\end{point}

As special case we get
\begin{point} Let $\{u(n)\}, \{v(n)\}$ be jointly naturally continuous sequences of elements of $[0,1]$ . Then
$$
u(n)=v(n), n \in \N
$$
If and only if $E(v)=E(u), E(u^2)=E(v^2), E(uv)=E(v^2)$.
\end{point}
\begin{point} If $\{v(n)\}, \{u(n)\}$ are sequences of elements from given interval $[a,b], a<b$, jointly naturally  continuous , having the same weak asymptotic distribution function and $E(uv)=E(v^2)$ then $u(n)=v(n)$ for $n \in \N$.
\end{point}

\section{Characterisation of naturality}

In the following we give a characterization when a natural measure
on given metric space exists. We say that a sequence $\{v(n)\}$ is
{\it strong positive } on set $S \subset \N$ if and only if there
exists a positive constant $k$ that $v(n)\ge k$ for $n \in S$.

\begin{point}
\label{natural} A metric $\mathfrak{d}$ is natural if and only if
for for every sequence uniformly continuous with respect the metric
$\mathfrak{d}$
\newline i) the limit $\lim_{N \to \infty}E_N(v)=E(v)$ exists and
\newline ii) $E(v)>0$ whenever $\{v(n)\}$ is non negative sequence
strongly positive on some closed ball with center from $\N$.
\end{point}

{\bf Proof.} Suppose that i) holds. If consider the space
$C(\MM)$ consisting from all continuous real valued functions
defined on $\MM$ then $E$ is a positive linear functional defined
on this space. Thus from Riesz representation theorem we get that
such probability Borelian measure $P$ defined on $\MM$ exists
that
$$
E(v)=\int \tilde{v} dP.
$$
This yields that $\{n\}$ is uniformly distributed with respect to
$P$.

Since the system of sets $B(n, \varepsilon), n \in \N,
\varepsilon>0$ forms an open base, it suffices to prove that $P(B(n,
\varepsilon))>0$ for given $n$ and $\varepsilon$. If $\delta <
\varepsilon$ then $B[n, \delta] \subset B(n, \varepsilon)$. Thus
Uryshon's theorem provides that such continuous function $f: \MM
\to [0,1]$ exists that $f(\alpha)=1, \alpha \in B[n, \delta]$ and
$f(\alpha)=0, \alpha \in \MM \setminus B(n, \varepsilon)$. And so
$\mathcal{X}_{B(n, \varepsilon)}(\alpha) \ge f(\alpha), \alpha \in
\MM$. This implies
$$
P(B(n, \varepsilon))=\int _{B(n, \varepsilon)} dP\ge\int f
dP=E(f)>0.
$$
\qed

\begin{theorem}
\label{natural1}
A metric $\mathfrak{d}$ is natural if and only if \newline
i) for each open ball $S(n, \varepsilon):=B(n, \varepsilon)\cap \N$ such
$\varepsilon' \le \varepsilon$ exists that
$S(n, \varepsilon')\in \D$ and $d(S(n, \varepsilon'))>0$, and \newline
ii) for each $\varepsilon >0$ such system of disjoint sets
$S_1, \dots, S_k \in \D$ exists that $\bold{diam} \ S_i < \varepsilon$
and $\N =S_1 \cup \dots \cup S_k$.
\end{theorem}

{\bf Proof.} We shall apply \ref{natural}.
 Let us suppose that the conditions i) and ii) hold. If $\{v(n)\}$ is a sequence uniformly continuous with respect to $\mathfrak{d}$, then for given $\delta >0$ such $\varepsilon$ exists that
\begin{equation}
\mathfrak{d}(m, n) < \varepsilon \Rightarrow |v(n)-v(m)| < \delta, \ m,n \in \N.
\end{equation}
Let
$$
\N= S_1 \cup \dots \cup S_k
$$
be the decomposition provided by ii). Let us consider the sequence $\{v_\varepsilon(n)\}$ defined as
$$
v_\varepsilon(n)=\sum_{j=1}^k \mathcal{X}_{S_j}(n)v(n_j)
$$
for suitable $n_j \in S_j, j=1, \dots, k$. Clearly limit
$$
\lim_{N \to \infty}E_N(v_\epsilon)= \sum_{j=1}^k d(S_j)v(n_j)=E(v_\varepsilon)
$$
exists. For each $n \in \N$ such $j$ exists that
$n \in S_j$. So $\mathfrak{d}(n, n_j) < \varepsilon$. Clearly
$v_\varepsilon(n)=v(n_j)$. This yields
$|v(n)-v_\varepsilon(n)| < \delta$. Therefore
$$
|\limsup_{N \to \infty}E_N(v)-E(v_\varepsilon)| < \delta
$$
and
$$
|\liminf_{N \to \infty}E_N(v)-E(v_\varepsilon)| < \delta.
$$
And so
$$
|\limsup_{N \to \infty}E_N(v)- \liminf_{N \to \infty}E_N(v)| < 2\delta.
$$
For $\delta \to 0^+$ we get that the proper limit
$\lim_{N \to \infty}E_N(v)=E(v)$ exists. We proved that i) of
\ref{natural} is fulfilled.

Suppose that $\{v(n)\}$ is a non negative sequence uniformly continuous and for some $m \in S(m, r), r>0$ the inequality
$v(n)> \kappa >0$ holds. From i) we get that such $r'< r$ exists that $S(m,r') \in \D$ and $d(S(m,r')) >0$. Then
$$
E_N(v) \ge \frac{1}{N}\sum_{n< N, n \in S(m,r')} \kappa >0.
$$
Thus
$$
E(v)\ge \kappa d(S(m,r')).
$$

Suppose that $\mathfrak{d}$ is an natural metric. Denote
$\MM$ the completion of $\N$. Let $\varepsilon >0 $ is given, Suppose that $r$ is such numbder that $\frac{\varepsilon}{2} > r>0$ is given.
Since $\N$ is dense in $\M$ we have that
$\M=\cup_{n=1}^\infty B(n,r)$.  And so taking account that the metric space is compact for suitable $n_1, \dots, n_s$ we have
$$
\M=\bigcup_{i=1}^s B(n_i, r).
$$
From \ref{pcont} we have that for each $i=1,2, \dots, s$ such
$\frac{\varepsilon}{2} > r_i\ge r$ exists that $B(n_i, r_i)$ is the sets of $P$ - continuity. And so
$$
\M=\bigcup_{i=1}^s B(n_i, r_i).
$$
The property \ref{pcontalgebra} allows us step by step construct the sets of $P$ - continuity $B_1, \dots, B_k$
such that
$$
\M=\bigcup_{i=1}^k B_i
$$
such that $diam B_i < \varepsilon$. If we put
$S_i=\N \cap B_i$, $i=1, \dots, k$ then these sets fulfill
the condition ii).

Consider an open ball $S(n, \varepsilon)$. Then
$S(n, \varepsilon)=\N \cap B(n, \varepsilon)$. From
\ref{pcont} we get that such positve $\varepsilon' < \varepsilon$ exists that $B(n, \varepsilon')$ is a set
of $P$ - contunuity. Thus $S(n, \varepsilon') \in \D$ and
$d(S(n, \varepsilon'))=P(B(n, \varepsilon'))>0$. We see that i) holds. \qed

\section{Construction of natural metrics} We repeat the construction of metric from \cite{IPT}(see also \cite{PAS5}).

Let $\mathcal{E}_n=\{E_1^{(n)}, \dots, E_{k_n}^{(n)}\}$, $n=1,2,3,\dots$ be such system of decompositions of $\N$ that \newline
a) the sets $E_{j}^{(n)}$ have positive asymptotic density, \newline b) each set from $\mathcal{E}_n$ is union of sets from $\mathcal{E}_{n+1}$ \newline
c) and the intersection
$\cap_{n=1}^\infty E^{(n)}_{j_n}$ has at most one element.

Define the function $\psi_n(a,b)=0$ if $a,b$ belong to the same set of $\mathcal{E}_n$ and
$\psi_n(a,b)=1$ otherwise, for $a,b \in \N$. Put

\begin{equation}
\label{naturalmetric}
\mathfrak{d}(a,b)=\sum_{n=1}^\infty \frac{\psi_n(a,b)}{2^n}
\end{equation}
By the standard way can be verified
\begin{point}
$\mathfrak{d}$ is metric.
\end{point}
\begin{theorem}
\label{construction} $\mathfrak{d}$ is natural metric.
\end{theorem}
{\bf Proof.} The assertion follows directly from Theorem \ref{natural1} \qed
\begin{example} Put
$$
r+(m)=\{a \in \N; a \equiv r \pmod{m}\},
$$
where $m \in \N, r \in \Z$. If $E^{(n)}_j = j+(n!)$
for $n \in \N$, $j=0, \dots, n!-1$ then the metric given by
(\ref{naturalmetric}) is the polyadic metric, (see \cite{N}, \cite{N1}, \cite{PAS5}).
\end{example}

\begin{theorem} A bounded one to one sequence of real numbers from some interval $[a,b]$ dense in this interval is continuous with respect to a suitable
natural metric if and only if it has increasing weak asymptotic distribution function.
\end{theorem}

{\bf Proof.} Let $\{v(n)\}$ be a sequence of elements of given interval $[a,b]$. From \ref{wadf} we get that this sequence has increasing weak asymptotic distribution function .

Suppose that $\{v(n)\}$ has increasing weak asymptotic distribution function $F$. Let
$I^{(K}_j, j=1, \dots, k_K$   a system of intervals $I^{(K}_j =[a_j^{(K}, a_{j+1}^{(K})$, $ j=1,\dots, k_K-1$ and
$I^{(K)}_{k_K}=[a_{k_K-1}^{(K)}, b]$ where $a_1^{(K)}=a$ be a system of divisions of $[a,b]$ such that
$a_j^{(K)}, j=2, \dots, k_K-1$ are the points of continuity of $F$ for $K \in \N$. Suppose that
$|I^{(K)}_j| \to 0$ for $K \to \infty$ uniformly for $i=1, \dots, k_K$. If we denote
$E_j^{(K)}=v^{-1}(I^{(K)}_j)$. The metric $\mathfrak{d}$ given by (\ref{naturalmetric}) is an natural metric and $\{v(n)\}$ is continuous
with respect to $\mathfrak{d}$. \qed
\begin{example}
Let us consider the sequence $\{u(n)\}$, where $u(n)=\frac{1}{n}, n \in \N$. Then $D^2(u)=0$ and so from \ref{mis1} we deduce that this sequence is not naturally continuous. If $\{v(n)\}$ is a sequence uniformly distributed modulo $1$ then the sequence $\{v_1(n)\}$, where $v_1(n)=v(n)+\frac{1}{n}$, is uniformly distributed modulo $1$ also. Thus these sequences are naturally continuous but not jointly naturally continuous.
\end{example}

\section{Laws of big numbers}
The continuity of the asymptotic distribution function of the sequence $f(v)$ will play an important rule in this part. One sufficient condition is given in \ref{abscont}.  

We derive a second one. 
In the paper \cite{pascont} the following result is proven:
\begin{point}
\label{cadf}Let $S$ be a set dense in real line.
A sequence $\{v(n)\}$ have continuous asymptotic distribution function if and only if \newline
i) $v^{-1}((-\infty, x))$ belongs to $\D$ for $x \in S$ and \newline
ii) for each real number $a$ we have $\lim_{\eta \to 0^+} \overline{d}(v^{-1}((a-\eta, a+\eta)))=0$.
\end{point}
\begin{point} Let $f$ be such continuous real function  defined on interval $[a,b]$ that such division
$a=a_0 < a_1 < \dots < a_m=b$ exists that the function
$f$ is strictly monotone on $[a_i, a_{i+1}], i=0, \dots, m$.
If $\{v(n)\}$ is a function having continuous asymptotic distribution function, then the sequence $\{f(v(n)\}$ have continuous asymptotic distribution function also.
\end{point}

{\bf Proof.} We apply \ref{cadf}. Denote $f_i=f|_{[a_i, a_{i+1}]}$ - the restriction of $f$ on $[a_i,a_{i+1}]$. Then for each interval $I$ we have
$$
f^{-1}(I)= \bigcup_{i=0}^{m-1}f_i^{-1}(I).
$$
Thus
\begin{equation}
\label{shortinterval}
v^{-1}(f^{-1}(I))= \bigcup_{i=0}^{m-1}v^{-1}(f_i^{-1}(I)).
\end{equation}
From the monotonicity of $f_i$ we get that $f_i^{-1}(I)$ are bounded intervals and so $v^{-1}(f_i^{-1}(I))\in \D$, thus $v^{-1}(f^{-1}(I)) \in \D$. The condition i) of \ref{cadf} is fulfilled.

The function $f_i$ are uniformly continuous. From the monotonicity we get that $f_i^{-1}$ are uniformly continuous also. Thus
$v^{-1}(f_i^{-1}(I))$ converges to $0$ if $|I| \to 0^+$.
And so from (\ref{shortinterval} ) we get that ii) of \ref{cadf} if fulfilled also. \qed

Let $\{v_k(n)\}, k=1,2,3,\dots $ be a one to one sequences of
elements of interval $[a,b]$ having continuous increasing asymptotic distribution functions $F_k$.
\begin{point}
\label{naturalmetric1} The sequences $\{v_k(n)\}, k=1,2,3, \dots$ are jointly naturally continuous.
\end{point}
{\bf Proof.} We shall construct a system of decompositions as in the case of Theorem \ref{construction}, (this construction is used also in the paper \cite{pasclt}).

Let $\{v_1(n)\},
\dots, \{v_j(n)\}, \dots $ be statistical independent sequences of
elements of $[0,1]$ having continuous asymptotic distribution
function $F$.
 Denote
$$
I^{(m)}_i=\Big[\frac{i}{2^m},\frac{i+1}{2^m}\Big), i=1, \dots,
2^m-2, I^{(m)}_{2^m-1}=\Big[\frac{2^m-1}{2^m},1\Big].
$$
where $m=1,2,3, \dots$. Put
\begin{equation}
\label{eee} E^{(m)}_{i_1, \dots, i_m} = \bigcap_{k=1}^m
v_k^{-1}(I^{(m)}_{i_k}), \ 0\le i_1, \dots, i_m \le 2^m -1 .
\end{equation}
Clearly  every  set $E^{(m)}_{i_1, \dots, i_m}$ has asymptotic
density and
\begin{equation}
\label{eeee} d(E^{(m)}_{i_1, \dots, i_m}) = \prod_{j=1}^m
\Big(F_j\Big(\frac{i_j+1}{2^{m}}\Big)-F_j\Big(\frac{i_j}{2^{m}}
\Big)\Big)>0,
\end{equation}
for
$\ 0\le i_1, \dots, i_m \le 2^m -1$.

Let $\mathcal{E}_m$ be the system of all sets in the form
(\ref{eee}) for $m \in \N$. Theorem \ref{construction} provides that the metric given by (\ref{naturalmetric}) is natural. \qed

In the following text we suppose that $\mathfrak{d}$ is the metric constructed in the proof of \ref{naturalmetric}, $P$ - given natural measure and $\nu^\ast$ measure density defined by (\ref{measuredensity}).

From Theorem \ref{equivalenceI} and \ref{nindependence} we can conclude:
\begin{point} The random variables $\tilde{v}_1, \dots, \tilde{v}_k, \dots $ are independent.
\end{point}
The proof of this statement can be obtain by straightforward transcription of the proof of Theorem 7 of \cite{P-T1}. We can prove small generalization of Theorem 4 of \cite{P-T1}:

\begin{theorem} Let $f$ such continuous real valued function defined on $[a,b]$ that the sequences $\{f(v_j)\}, j=1, \dots, N$ have continuous asymptotic distribution function. Then the set
$$
V(N, \varepsilon)=
\Big\{n \in \N;\Big|\frac{\sum_{j=1}^N f(v_j(n))}{N}-\frac{\sum_{j=1}^N \int_a^b f(x)dF_j(x)}{N}\Big| \le \varepsilon \Big\}
$$
belongs to $\D_\nu$ and
$$
d(V(N,\varepsilon))=\nu(V(N,\varepsilon))\ge 1- \frac{C^2}{N\varepsilon^2}
$$
for each $\varepsilon >0$, where $C=\sup\{|f(x); x \in [a,b]\}$.
\end{theorem}

{\bf Proof.} We have $E(f(v_j))=\int_a^b f(x)dF_j(x)$, where
$j=1, \dots, N$. The sequences $\{f(v_j)\}$ have continuous asymptotic distribution function. This yields that
the sequence  $\{\sum_{j=1}^N f(v_j(n))\}$ has continuous asymptotic distribution function. Moreover this sequence is uniformly continuous with respect to $\mathfrak{d}$. The random variables $f(\tilde{v}_j), j=1,\dots, N$ are independent and $E(f(\tilde{v}_j)= E(f(v_j))$. Moreover
$D^2(\tilde{v}_j)\le C^2$.
Thus aplying \ref{densint} and weak law of big numbers we get the assertion. \qed

\begin{point}
\label{slbn}If $f$ is such real continuous function defined on $[a,b]$ that the sequences $\{f(v_i(n)\}$ have continuous asymptotic distribution function, then
for almost all $\alpha \in \MM$ we have
$$
\lim_{N \to \infty} \frac{f(\tilde{v}_1(\alpha))+\dots + f(\tilde{v}_N(\alpha))}{N}- \frac{\sum_{j=1}^N\int_a^bf(x)dF_j(x)}{N}=0.
$$
\end{point}

\begin{point}Suppose that $f$ is real valued continuous monotone function defined on $[a,b]$. For almost all
$\alpha \in \MM$ exists non increasing sequence of $\varepsilon_n >0 ; n \in \N$ such for each sequence
of natural numbers $\{k_n\}$ where $k_n \in B(\alpha, \varepsilon_n)$ we have
$$
\lim_{N \to \infty} \frac{f(v_1(k_1))+\dots + f(v_N(k_N))}{N}- \frac{\sum_{j=1}^N\int_a^bf(x)dF_j(x)}{N}=0.
$$
\end{point}

\begin{theorem} If the sequences $\{v_i(n)\}, i=1,2,3, \dots$  are uniformly distributed modulo $1$ then the sequence
$\{v_\alpha(n)\}$ where $v_\alpha(n)=\tilde{v}_n(\alpha), n\in \N$
is uniformly distributed modulo $1$ for almost all $\alpha \in \M$.
\end{theorem}

{\bf Proof.} Denote
$$
S_h =\Big\{\alpha \in \M; \lim_{N \to \infty} E_N(v^h_\alpha)=
\frac{1}{h+1} \Big\}
$$
for $h =1,2,3, \dots$. The criterion \ref{wkpolynom} yields that the sequence $\{v_\alpha(n)\}$ is uniformly distributed modulo $1$ if and only if $\alpha \in
S_h$ for all $h=1,2,3, \dots $. From \ref{slbn} we get
$P(S_h)=1$. Thus $P(\cap_{h=1}^\infty S_h)=1$. \qed

\begin{point} For almost all $\alpha \in \M$ such decreasing
sequence $\varepsilon_n >0$ exists that for each sequence of natural numbers $\{k_n\}$ such that $k_n \in B(\alpha,\varepsilon_n)$ the sequence $\{v_n(k_n)\}$ is uniformly distributed modulo $1$.
\end{point}

{\bf Proof.} The functions $\tilde{v}_n, n=1,2,3, \dots $ are continuous. Therefore such sequence $\varepsilon_n >0$ exists that $\mathfrak{d}(\alpha, a) < \varepsilon_n \Rightarrow
|v_n(m)-\tilde{v}_n(\alpha)| < \frac{1}{2^n}$. Thus if
$k_n \in B(\alpha, \varepsilon_n)$ then
$\lim_{n \to \infty} v_n(k_n)-\tilde{v}_n(\alpha) =0$. \qed

Let us remark that the lasts proofs are straight forward generalization of the polyadic cases from \cite{P-T1}.

Author's adress: Department of Mathematics and Informatics, Faculty of Education, University of Trnava,
Priemyselna 4, Trnava, Slovakia.

\end{document}